\documentclass[12pt]{amsart}

\newtheorem{theorem}{Theorem}[section]
\newtheorem{proposition}[theorem]{Proposition}
\newtheorem{corollary}[theorem]{Corollary}
\newtheorem{remark}[theorem]{Remark}
\newtheorem{lemma}[theorem]{Lemma}

\usepackage{amssymb}

\numberwithin{equation}{section}

\begin{document}

\baselineskip=16pt

\title[coherent systems]{Tensor product of coherent systems}

\author{L. Brambila-Paz}

\address{CIMAT, Apdo. Postal 402, C.P. 36240.
 Guanajuato, Gto,
M\'exico}

\email{lebp@cimat.mx}

\author{Angela Ortega }

\address{Instituto de Fisico-Matematicas, Universidad Michoacana, Morelia, Mich. , M\'{e}xico}

\email{ortega@ifm.umich.mx}

 \keywords{coherent systems, moduli space,
stability, Brill-Noether, tensor product}

\subjclass[2000]{14H60, 14J60}

\date{}
\begin{abstract}
Let $X$ be a smooth algebraic curve of genus $g\geq 2$. A stable
vector bundle over $X$ of degree $d$, rank $n$ with at least $k$
sections is called a Brill-Noether bundle of type $(n,d,k)$. By
tensoring coherent systems, we prove that most of the known
Brill-Noether bundles define coherent systems of type $(n,d,k)$
that are $\alpha$-stables for all allowable $\alpha .$
\end{abstract}

\maketitle

\section{Introduction}

Let $X$ be a smooth projective algebraic curve over $\mathbb{C}$
of genus $g\geq 2.$ By a coherent system of type $(n,d,k)$ we mean
a vector bundle of degree $d$ and rank $n$ together with a linear
subspace of its space of sections of dimension $k$.

 In \cite{lep}, \cite{an}, \cite{rag} was
introduced a notion of stability for coherent systems which
permits the construction of moduli spaces. This notion depends on
a real parameter $\alpha$, and thus leads to a family of moduli
spaces. The different members in the family correspond to
different values of $\alpha .$ As $\alpha $ varies, the moduli
spaces can change only when $\alpha$ passes through one of a
discrete set of points, called critical values. If $k<n$, the
range of parameter is a finite interval and the family has only
finitely many distinct members. If $k\geq n$ the range of the
parameter is infinite; however, there is only a finite number of
distinct moduli spaces. Moreover, if $G(\alpha :n,d,k)$ is the
moduli space of $\alpha$-stable coherent systems of type
$(n,d,k)$, there is a critical value, denoted by $\alpha _L$, such
that for all $\alpha , \alpha '
>\alpha _L$, $G(\alpha :n,d,k)=G(\alpha ' :n,d,k)$
(\cite{bomn}).

Let $B(n,d,k)$ be the Brill-Noether locus of the stable vector
bundles over $X$ of degree $d$, rank $n$ with at least $ k$
independent sections. A triple $(n,d,k)$ is called a Brill-Noether
triple if $B(n,d,k)\not= \emptyset $ and a vector bundle in
$B(n,d,k)$ is called a Brill-Noether bundle of type $(n,d,k)$.  It
is clear that Brill-Noether bundles define coherent systems and
for $\alpha
>0$ close to $0,$ they define $\alpha$-stable coherent systems
(see \cite{bo}). However, for arbitrary choice of $\alpha$ it is
not clear the relation between Brill-Noether bundles and stable
coherent systems. We are interested in studying this relationship.

The case $k\leq n$ have already been considered in \cite{bo}
\cite{bomn} and \cite{bommn} and was proved that if $n\leq d
+(n-k)g$ and $(n,d,k)\not= (n,n,n),$ Brill-Noether bundles define
 $\alpha$-stable coherent systems for all allowable $\alpha .$
Moreover, the set of such bundles is a Zariski open subset of
$G(\alpha :n,d,k).$

For general curves, if  $k=n+1$ and $d\leq g+n$ it was proved in
\cite{lb} that there exist Brill-Noether bundles which are $\alpha
$-stable for all $\alpha >0$. The existence of such bundles was
used to determine the structure of $G(\alpha :n,d,k)$ and to prove
that $\alpha _L =0$. For $k>n$ and $n,d,k$ satisfying the
relations in \cite{mon1}, Montserrat Teixidor i Bigas in
\cite{mon} constructs a singular curve and a limit of coherent
system such that, on a generic curve, define a coherent system
that is $\alpha$-stable for all $\alpha
>0$ and the vector bundle is stable.

Our aim in this paper is to study the $\alpha$-stability of
coherent systems defined by the Brill-Noether bundles given in
\cite{me} and \cite{bfmn}. Such bundles will be used to determine
the structure of $G(\alpha :n,d,k)$ and, in some cases, to
determine the value of $\alpha _L$. Our results for $k>n$ will be
for any curve and, even for the generic case, they extend beyond
of those in \cite{mon}. Furthermore, our methods give another
proof of non-emptiness for those values that are included.

In order to state our results we recall the following definitions
and facts.

Denote by $\mathcal{M}(n,d) $ the moduli space of stable vector
bundles of rank $n$ and degree $d$. If $d>n(2g-2)$, denote by
$Grass(s)$ the Grassmannian bundle over $\mathcal{M}(n,d)$ with
fibre the Grassmannian $Grass (s,H^0(E))$.

From \cite[Corollaire 3.14]{he}, every irreducible component of
$G(\alpha :n,d,k)$ has dimension greater or equal to the
Brill-Noether number $\beta(n,d,k):=n^2(g-1)+1-k(k-d+n(g-1)).$
Denote by $G_0(n,d,k)$ the first member of the moduli spaces
family and by $G_L(n,d,k)$ the last one. For coherent systems
$(E,V)$ of type $(n,d,k)$ with $k\geq 1$ define $G_g(n,d,k)$ and
$U(n,d,k)$
 as
 $$G_g(n,d,k):=\{(E,V): (E,V)
 \ \mbox{is generated by $V$ with $H^0(E^*)=0$}\}.$$
 and
$$
U(n,d,k):=\{(E,V):(E,V) \, \, \mbox{is  }  \, \alpha -\mbox{stable
for all allowable} \, \alpha \mbox{ and }E\mbox{ is stable}\},
$$

We are interested in studying the non-emptiness of $U(n,d,k),$
when $(n,d,k)$ is a Brill-Noether triple.

We now state our results.

Let  $2 \leq gs <n$ and $k=n+s$. If $d<n+gs$, it is known (see
\cite{me}) that there are no semistable bundles of degree $d$,
rank $n$ with at least $k$ independent sections; hence
$G_0(n,d,k)=\emptyset$. Moreover, the non existence of semistable
bundles also implies that (see Theorem \ref{teoempty})

\begin{itemize}\item {\it if $d<n+gs$, $G(\alpha:n,d,k)= \emptyset $ for all $\alpha
>0$.}
\end{itemize}
If $d\geq n+gs$, from \cite{me} there are Brill-Noether bundles of
type $(n,d,k)$. In this case we prove

\begin{itemize}
\item {\it If  $d=n+gs$, {\mbox{(see Theorem \ref{teo1})}}}
\end{itemize}
 \begin{enumerate} {\it
 \item $G(\alpha:n,d,k)\not= \emptyset$ for all
$\alpha >0$;
 \item $G(n,d,k): = G(\alpha: n, d, k) = G(\alpha':n,
d, k) $ for $\alpha ,\alpha ' >0$ i.e. $\alpha _L =0$;\item
$U(n,d,k)\not= \emptyset $. Moreover, $U(n,d,k)=G(n,d,k)\subset
G_g(n,d,k);$ \item $ G(n,d,k)$ is smooth irreducible of dimension
$ \beta (n,d,k) $. Moreover, $G(n,d,k)\cong \mathcal{M}(s,d) $.
\item For $s\geq 1$, $G(\alpha:s,d,k)= G(n,d,k)$ for all $\alpha
>0$.}
\end{enumerate}

\begin{itemize}\item
{\it For $d=n+sg+s'$ with $0<s'<g$ {\mbox{(see Theorem
\ref{teo2})}};}
\end{itemize}

\begin{enumerate} {\it \item $G(\alpha:n,d,k)\not= \emptyset$ for all
$\alpha >0$; \item $U(n,d,k)\not= \emptyset$. Moreover, any
$(E,V)\in U(n,d,n+s)$ is generically generated; \item $ G(\alpha
:n,d,k)$ has a smooth irreducible component $G^0(n,d,k)$ of
dimension $ \beta $. Moreover, $G^0(n,d,k)$ is birationally
equivalent to the Grassmanian bundle $Grass (s')$ over
$\mathcal{M}(s,d)$.}
\end{enumerate}

From Theorem \ref{teo1}, \ref{teo2} and \cite[Theorem A]{bommn} we
have that
\begin{itemize} \item {\it If $0<d<2n$, there is an open set
$Z$ of  $B(n,d,k)$ such that any $E\in Z$ defines a coherent
system that is $\alpha$-stable for all $\alpha
>0$. Moreover, for any Brill-Noether triple $(n,d,k)$,
$U(n,d,k)\not= \emptyset .$} \end{itemize}

 The idea of tensoring Brill-Noether
bundles by line bundles with section was used in \cite{bfmn} to
produce Brill-Noether bundles of degree $d>2n.$ We use this idea
and tensor coherent systems of type $(n,d,k)$ by coherent systems
of type $(1,d',k')$.

 Let $(E,V)$ be a coherent system  and $L$ an
effective line bundle. Choose a section $s$ of $L$ and define the
coherent system $(E\otimes L,\widetilde{V})$ where $\widetilde{V}$
is the image of $V$ in $H^0(E\otimes L)$ under the canonical
inclusion $H^0(E) \hookrightarrow H^0(E\otimes L)$ induced by $s$.
Raghavendra and Vishwanath in \cite[Lemma 1.5]{rag} proved that
$(E,V)$ is $\alpha $-stable if and only if $(E\otimes L,
\widetilde{V})$ is $\alpha $-stable. We extend such Lemma as
follows.

Let $(L,W)$ be a coherent system of type $(1,d',k')$ and $K_L$ the
kernel of the evaluation map $W\otimes \mathcal{O} \rightarrow L$.
Let $(E,V)$ be a coherent system of type $(n,d,k)$. If
$H^0(K_L\otimes E)=0$, we identify $W\otimes V$ with the image of
$W\otimes V$ in $H^0(E\otimes L)$ under the
 inclusions $W\otimes V\hookrightarrow W\otimes H^0(E)
\hookrightarrow H^0(E\otimes L)$. Hence, $(E\otimes L, V\otimes
W)$ is a coherent system of type $(n,d+nd',kk')$. For such
coherent systems we prove (see Lemma \ref{proplemma})
\begin{equation}\label{eq1}\end{equation}\begin{itemize} \item
{\it $(E,V)$ is $k'\alpha$-stable if and only if $(E\otimes L,
V\otimes W)$ is $\alpha$-stable.}\end{itemize}

We use $(\ref{eq1})$ and Theorems \ref{teo1}, \ref{teo2} and
\cite[Theorem A]{bommn} to prove that most of the Brill-Noether
bundles in \cite{bfmn} are $\alpha $-stables for all allowable
$\alpha.$ Actually, we prove that

{\it under the hypothesis that $d=d''+nd'$ and $k=k'k''$ with
$0<d''<2n$, $n\leq d''+(n-k'')g$ and $(n,d'',k'')\not= (n,n,n)$;
 and $(d',k')$ satisfying one of the following conditions }
\begin{enumerate}
\item {\it $d'\leq  2g$ and $\beta (1,d',k')\geq 0$ (see Theorem
\ref{teogrande02}); \item $d' > 2g$ and $k'\geq 1$ and
$nd'>(k-1)d''$ (see Theorem \ref{teodetodos})},\end{enumerate}

\begin{itemize}
\item {\it $G(\alpha:n,d,k)\not= \emptyset $ for all allowable
$\alpha $. Moreover, $U(n,d,k)\not= \emptyset$}\end{itemize}

For a hyperelliptic curve, under the same hypothesis on
$(n,d'',k'')$ and the assumptions that $k'\geq 1$ and $d'=2(k'-1)$
we prove
\begin{itemize}
\item {\it $G(\alpha:n,d,k)\not= \emptyset $ for all allowable
$\alpha $. Moreover, $U(n,d,k)\not= \emptyset$}.\end{itemize}

During the final stages of writing this paper, we came across the
work in \cite{bommn2}. Our results in Theorem \ref{teo1} and
\ref{teo2} partly coincide with some of their results.
\smallskip

{\it Notation}

We will denote by $K$ the canonical bundle over $X$,
 by $K_E$ and $I_E$ the kernel and image, respectively,
of the evaluation map $V\otimes \mathcal{O}\rightarrow E$,
 $H^i(X,E)$ by $H^i(E)$, $\dim H^i(X,E)$
by $h^i(E)$, the rank of $E$ by $n_E$, the degree of $E$ by $d_E$.
If $W$ is a vector space, we denote $W\otimes \mathcal{O}\otimes
E$ by $W\otimes E$. If $\beta: V\rightarrow H$ and $\gamma :
W\rightarrow H$ are injective homomorphism, we will denote $\beta
(V) \cap \gamma (W)\subset H$
 by $V\cap W$.

\bigskip

\section{Definitions and general results}

We recall some definitions and facts of coherent systems we shall
need. We refer the reader to \cite{bomn} and \cite{bdo} and
 references cited therein for basic properties of coherent
systems on algebraic curves.

Let $X$ be a smooth projective algebraic curve over $\mathbb{C}$
of genus $g\geq 2$. A coherent system over $X$ of type $(n,d,k)$
is a pair $(E,V)$ where $E$ is a vector bundle over $X$ of rank
$n$, degree $d$ and $V$ a linear subspace of $H^0(X,E)$ of
dimension $k$. For any real number $\alpha > 0,$ define the
$\alpha$-slope of the coherent system $(E,V)$ of type $(n,d,k)$ as
$$\mu _{\alpha}(E,V):= \mu(E) + \alpha \frac{k}{n},$$
where $\mu (E):=d/n$ is the slope of the vector bundle $E.$ A
coherent subsystem $(F,W)\subseteq (E,V)$ is a
 coherent system such that $F\subseteq E $ and $ W\subseteq V\cap
 H^0(F).$
For any $\alpha $, a coherent system $(E,V)$ is $\alpha$-stable
(respectively $\alpha$-semistable) if for all proper coherent
subsystems $(F,W)$
$$\mu _{\alpha}(F,W)<\mu _{\alpha}(E,V) \ \ \
({\rm respectively} \leq ).$$

 Denote by $G(\alpha:n,d,k)$
(respectively $\widetilde{G}(\alpha:n,d,k)$) the moduli space of
$\alpha $-stable (respectively $\alpha$-semistable) coherent
systems of type $(n,d,k)$. For non-emptiness of $G(\alpha:n,d,k)$
with $k\geq 1$ we need $\alpha
>0$ and $d>0$. Basic properties of $G(\alpha:n,d,k)$ have been proved in
\cite{lep}, \cite{an} and \cite{rag}.

Most of the detailed results known are for $k\leq n$ (see
\cite{bo}, \cite{bommn},\cite{bomn}). For $k=n+1$ and $X$ general
see \cite{lb}, \cite{bun} and \cite{bomn}. For $k>n$, on a generic
curve, Montserrat Teixidor i Bigas in \cite{mon} proved that,
under the same relation as in \cite{mon1}, $G(\alpha:n,d,k)\not=
\emptyset$ and has an irreducible component of the correct
dimension. For $d>>0$ see \cite{ba}. It is our purpose here to
study the case $k>n$ on any curve and under different conditions
then those in \cite{mon}.

Every irreducible component of $G(\alpha:n,d,k)$ has dimension at
least the Brill-Noether number $\beta (n,d,k):= n^2(g-1)+1
-k(k-d+n(g-1))$. From the infinitesimal study of the coherent
systems if $(E,V)\in G(\alpha:n,d,k)$, $G(\alpha:n,d,k)$
 is smooth of dimension $\beta (n,d,k) $
 in a neighbourhood of $(E,V)$ if
and only if the {\it Petri map} $$V\otimes H^0(E^*\otimes
K)\rightarrow H^0(E\otimes E^*\otimes K)$$ is injective.

Given a triple $(n,d,k)$ denote by $C(n,d,k)$ the set
$$C(n,d,k):=\{\alpha | 0\leq \alpha
=\frac{nd'-n'd}{n'k-nk'} \  {\rm with} \ \ 0\leq k'\leq k,
0<n'\leq n, \ {\rm and} \ \
 nk'\not= n'k \}.$$

An element $\alpha $ in $ C(n,d,k)$ is
 called a critical point. The set $C(n,d,k)$  defines a
 partition of the interval $[0,\infty )$. With
 the natural order on $\mathbb{R}$, label the critical points as
$\alpha _i.$

It is known (see \cite{bdo} and \cite{bomn}) that

\begin{enumerate}
\item If $\alpha ',\alpha '' \in (\alpha _i , \alpha _{i+1})$ then
  $G(\alpha ':n,d,k)= G(\alpha '':n,d,k)$.
Denote by
 $G_i(n,d,k)$ the moduli space $G(\alpha:n,d,k)$
for any $\alpha \in (\alpha _i , \alpha _{i+1}).$ \item For $k<
n$, if $\alpha > \frac{d}{n-k}$, $G(\alpha:n,d,k)=\emptyset$.
\item For $k\geq n,$ there exists $\alpha _L$ such that for any
$\alpha , \alpha '>\alpha _L, \ \ G(\alpha:n,d,k)= G(\alpha
':n,d,k)$. Denote by $G_L(n,d,k)$ the moduli space
 $G(\alpha:n,d,k)$ for $\alpha > \alpha _L$.
\end{enumerate}

Let $B(n,d,k)$ (respectively $\widetilde{B}(n,d,k)$) be the
Brill-Noether locus of stable (respectively semistable) vector
bundles. There is a natural map
$$\phi :  G_0(n,d,k) \rightarrow  \widetilde{B}(n,d,k)$$
defined by $(E,V) \mapsto E$ that is injective over
$B(n,d,k)-B(n,d,k+1).$ Moreover, if $E\in B(n,d,k)$ then for any
subspace $V\subseteq H^0(E)$ of dimension $k$, $(E,V)\in
G_0(n,d,k)$.

\begin{remark}\begin{em}\label{rem1} Let $(E,V)$ be a coherent system of type
$(n,d,k)$. From the definition of $\alpha $-stability and
stability of a vector bundle we have that
\begin{enumerate}
\item if $n=1$, $(E,V)$ is $\alpha $-stable for all $\alpha >0$;
\item if $(E,V)\in G(\alpha:n,d,k)$ and $E$ is stable then $(E,V)$
is $\alpha '$-stable for all $0<\alpha ' <\alpha ;$ \item if $E$
is stable and for all subsystems $(F,W) \subset (E,V)$,
$\frac{{\rm dim}W}{n_F} \leq \frac{k}{n}$ then $(E,V)$ is
$\alpha$-stable for all $\alpha >0.$
\end{enumerate}\end{em}
\end{remark}

 For coherent systems of type $(n,d,k)$ define $U(n,d,k)$ as
$$
U(n,d,k):= \{ (E,V): (E,V)\  \ {\mbox is}
 \, \  \alpha {\mbox {-stable  for  all allowable $\alpha$  and E {\rm is stable}}} \ \
 \}.$$

By \`{}\`{}allowable\'{}\'{} $\alpha $ we mean that if $k<n$,
$\alpha <\frac{d}{n-k}$ and if $k\geq n$, $\alpha >0.$

Note that $G(\alpha :n,d,k)$ could be non-empty for all $\alpha
>0$, but $U(n,d,k)=\emptyset.$ From the openness of
$\alpha$-stability $U(n,d,k)$ is an open subset of $G_L(n,d,k)$.
Moreover, if $G_L(n,d,k)$ is irreducible, $U(n,d,k)$ is
irreducible.

If $k\leq n$ and $n\geq 2$, $U(n,d,k)\not= \emptyset$ if and only
if $n\leq d +(n-k)g$ and $(n,d,k)\not= (n,n,n)$ (see \cite[Theorem
A]{bommn}). Hence, the Brill-Noether bundles with $k\leq n$ are
$\alpha $-stable for all allowable $\alpha .$

Our aim is to study $U(n,d,k)$ for $k>n$. In particular, the
non-emptiness.

A coherent system $(E,V)$ can be defined as a triple
$(E,V,\phi_{E,V})$ where $V$ is a vector space and $\varphi
_{E,V}: V\otimes \mathcal{O}\rightarrow E$ is a map such that the
induced map $\varphi ^*_{E,V}: V\rightarrow H^0(E)$ is injective.
Moreover, we have the exact sequence
\begin{equation}\label{eqsistems}
0\rightarrow K_E \rightarrow V\otimes \mathcal{O}
\stackrel{\varphi _{(E,V)}}{\rightarrow} E\rightarrow H\oplus
\tau\rightarrow 0
\end{equation}
where $K_E$ and $H$ are vector bundles with $H^0(K_E)=0$ and
$\tau$ a torsion sheaf.

If $I_E$ is the image of the evaluation map $V\otimes
\mathcal{O}\rightarrow E$ we split $(\ref{eqsistems})$ as
\begin{equation}\label{eqsystems1}
0\rightarrow K_E \rightarrow V\otimes \mathcal{O}\rightarrow I_E
\rightarrow 0
\end{equation}
and
\begin{equation}\label{eqsystems2}
0\rightarrow I_E \rightarrow E\rightarrow H\oplus \tau \rightarrow
0.
\end{equation}

If $K_E=0$, the coherent system is called injective. If $H=0$, is
called generically generated and if also $\tau=0$, generated. Note
that if $(E,V)$ is generically generated, the rank of $I_E$ is
$n.$
\begin{remark}\begin{em}
If $(E,V)$ is generated, the dual of the kernel of the evaluation
map (i.e. the vector bundle $K_E^*$ in $(\ref{eqsistems})$) is
usually  denoted as $M_{V,E}$ and if $V=H^0(E)$, as $M_E$.
\end{em}\end{remark}

\begin{remark}\begin{em}\label{remlineal}
If $k \leq n$, recall from \cite [Corollary 2.5]{bommn} that every
$\alpha $-semistable coherent system is injective, if $d\leq min
\{2n,n+\frac{ng}{k-1}\}.$ Moreover, for any injective coherent
system $(E,V)$, $H^0(H^*)=0$. In particular, if $k=n$, $H=0.$
\end{em}\end{remark}

 If $k> n$, from \cite[Proposition 4.4]{bomn},
there exists $\alpha _{gg}$ such that for any $\alpha > \alpha
_{gg}$, a $\alpha $-semistable coherent system is generically
generated. Actually, $\alpha _{gg} \leq \frac{d(n-1)}{k}.$

\begin{remark}\begin{em}\label{remdeu} Note that if $k>n$ and
$(E,V)\in U(n,d,k)$, $E$ is semistable and generically generated.
\end{em}\end{remark}

For any $(n,d,k)$ define $G_g(n,d,k)$ as
$$G_g(n,d,k):=\{(E,V): (E,V) \mbox{ is of type } (n,d,k) \, \mbox{ and it is
generated with $H^0(E^*)=0$}\}$$

\begin{remark}\begin{em}\label{rem11} Let $(E,V)\in G_g(n,d,k)$. From
\cite[Proposition 2.5]{lb},
 {\begin{enumerate} \item if
$G_g(n,d,k)\not=\emptyset$, $k>n;$ \item any quotient coherent
system $(Q,Z)$ is generated; \item  if $k=n+s$ with $s\geq 1$, for
any subsystem $(F,W)\subset (E,V)$, $h^0(F)\leq n_F+s-1$; \item if
$d\geq 2gn$, $G_g(n,d,k)\not=\emptyset$ for $k\geq n+1$;\item  if
$E$ semistable and $k=n+1$, $(E,V)$ is $\alpha$-stable for all
$\alpha
>0.$\end{enumerate}}
\end{em}
\end{remark}

Given $(G,V)\in G_g(n,d,k)$, we have the exact sequence
\begin{equation}\label{eqdual}
0\rightarrow M_{V,G}^* \rightarrow V\otimes \mathcal{O} \rightarrow
G\rightarrow 0
\end{equation}
with $H^0(M^*_{V,G})=0.$ The vector bundle $M_{V,G}$ is generated
by $V^*$ and $(M_{V,G},V^*)$ is called the {\it dual span} of
$(G,V)$. Moreover, if $(n,k)=1$, $(G,V)$ is $\alpha$-stable for
large $\alpha $ if and only if $(M_{V,G},V^*)$ is $\alpha$-stable
for large $\alpha $ (see \cite[Corollary 5.10]{bomn}).

\begin{remark}\begin{em}\label{remlineal1} In particular, if $(L,V)\in G_g(1,d,k)$,
the dual span $(M_{V,L},V^*)$ is $\alpha$-stable for large $\alpha
$, since $(L,V)$ is $\alpha $-stable for all $\alpha$.
\end{em}\end{remark}

Let $G$ be a stable of degree $d>2gn$. In \cite{bu1} Butler proved
that $M_{G} $ is stable and Mercat in \cite{me} gives an
isomorphism between $G_0(n,n+sg, n+s)$ and $G_0(s,n+sg, n+s)$ for
$sg<n$. In section \S 4 we will study this relation for all
$\alpha >0.$

If $X$ is general and $g>n^2-1,$ using the dual span
correspondence, Butler in \cite{bu} gave a birrational map between
$G_0(n,d, n+1)$ and $G_0(1,d, n+1)$. The dual span correspondence
was also used in \cite[Theorem 5.11]{bomn} to give necessary and
sufficient conditions for non-emptiness of $G_L(n,d,n+1)$ and in
\cite[Theorem 4.7]{lb} to prove non-emptiness and to describe the
structure of $G(\alpha:n,d,n+1)$ for all $\alpha >0.$

\begin{remark}\begin{em}\label{remdualbundle} The vector bundles
$M_{V,G}$ have been studied from different points of view (see
e.g. \cite{le}, \cite{pr} \cite{beu}). The existence of line
bundle of degree $d$ with $M_{V,L}$ stable and $\dim V=n+1$, has
been proved in the following cases;
\begin{enumerate}
\item $d\geq 2g$ and $ d+(1-g)\geq n+1\geq 2$ (see \cite{me}and
\cite{bu1}).\item If $K$ is the canonical bundle and $X$ is
non-hyperelliptic (see \cite{pr}).\item If $X$ is general and $n+g
-\frac{g}{n+1}\leq d\leq g+n$ (see \cite{lb}, \cite{bu} and
\cite{ye}). \item If $X$ is general and; $d\geq n+g
-\frac{g}{n+1}$ and $g\geq n^2-1$ (see \cite{lb}and \cite{bu}).
\end{enumerate}
\end{em}\end{remark}

\begin{remark}\begin{em}\label{rem111} Note that if
$G_g(n,d,k)\not=\emptyset ,$ from \cite[Proposition 3.2]{pr},
$G_g(n,d,n+1)\not=\emptyset $ and by the dual span correspondence
$G_g(1,d,n+1)\not=\emptyset .$
\end{em}
\end{remark}

\section {Brill-Noether Bundles}

In this section, we recall some facts and the construction of the
Brill-Noether bundles in \cite{bgn}, \cite{me} and \cite{bfmn}.

First we recall from \cite{me} a Proposition that will be used.

\begin{proposition}\label{propm} \cite[Proposition A.2]{me}
Let $F$ be a vector bundle of rank $n$ and degree $d$ such that
$H^0(F^*)=0$. If the maximal semistable subbundle of $F$ has slope
$<2$ then $h^0(F)\leq n+\frac{d-n}{g}.$
\end{proposition}

 Let $E$ be a stable bundle of degree $d$ and rank $n\geq 2$.
\begin{itemize}
\item  If $d>n(2g-2),$ $U(n,d,k)=\emptyset$ for $k>d+n(1-g)$;
\item if $d\geq 2gn$, $G_g(n,d,k)\not= \emptyset;$ \item if
$0<d<n$ and $k>n$, $U(n,d,k)=\emptyset$ (see \cite{bgn}).
\end{itemize}

\begin{remark}\begin{em}\label{remstabmerc}
If $0<d< 2n$, there exist stable vector bundles of rank $n$ and
degree $d$ with $k$ independent sections if and only if $n\leq
d+(n-k)g$ and $(n,d,k)\not=(n,n,n)$ (see \cite{bgn} and \cite{me}).
Therefore, $G_0(n,d,k)\not= \emptyset $ if and only if $n\leq
d+(n-k)g$ and $(n,d,k)\not=(n,n,n)$. Moreover, if $n>d+(n-k)g$,
$U(n,d,k)=\emptyset .$
\end{em}\end{remark}

If $0<d<2n$, stable bundles with $k\leq n$ define injective
coherent systems that are $\alpha $-stable for all allowable
$\alpha$ (see \cite{bomn}, \cite{bo} and Remark \ref{remlineal}).
Hence $U(n,d,k)\not=\emptyset .$

For $k>n$,  we know from \cite[2-B1]{me} that any such stable
bundle $A$ fits in an exact sequence
\begin{equation}\label{eqa}
0\rightarrow G^* \rightarrow H^0(A)\otimes \mathcal{O}\rightarrow
A\rightarrow 0.
\end{equation}
where $G$ is a stable bundle of rank $s$, slope $>2g$ and from
Proposition \ref{propm}, $h^0(G)=h^0(A)=k.$ Actually, $A=M_{G}$
and the coherent system $(A,H^0(A))$ is in $G_0(n,d,k)$ and it is
generated. Moreover, $(G,H^0(G))\in G_0(s,d,k)$.

If $W$ is a general subspace of $H^0(A)$ of $\dim W=s'$, we have
the exact sequence
\begin{equation}\label{eqkernel2}
0\rightarrow W\otimes \mathcal{O} \rightarrow A\rightarrow
B\rightarrow 0.
\end{equation}

From \cite[3-B1]{me} and its proof, we know that any such bundle
$B$ that fitting  in the exact sequence $(\ref{eqkernel2})$ is
stable.

\begin{remark}\begin{em}\label{remcondit} Note that the condition $d>2gs$
with $0<d<2n$ is equivalent to $0<gs<n$ and $d\geq n+gs.$
Moreover, if $k=n+s$,
$$d\geq n+gs \, \, \Longleftrightarrow \, \, n\leq d+(n-k)g.$$
\end{em}\end{remark}

\begin{remark}\begin{em}\label{remdgrande}
In particular, if $d=(n+s')+sg$, with $0<s'<g$,
$h^0(G)=h^0(A)=(n+s') +s$. Moreover, $n_B=n$, and from Proposition
\ref{propm}, $h^0(B)=n+s$ and $H^0(A)=W\oplus H^0(B)$. The
coherent system $(B,H^0(B))$ is in $G_0(n,d,k)$ and it is
generated.
\end{em}\end{remark}

From the cohomology sequence of $(\ref{eqkernel2})$ and Remark
\ref{remdgrande} we have that the coboundary map $\delta : H^0(B)
\rightarrow H^1(\mathcal{O})\otimes W$ is the zero map and hence,
we have the following diagram

\begin{equation}\label{diagram1}
\begin{array}{ccccccccc}
&& && 0&
& 0&&\\
&& && \downarrow&
& \downarrow&&\\
&& && G^{{*}}&=
& G^{{*}}&&\\
&& && \downarrow&
& \downarrow&&\\
0&\rightarrow& W\otimes \mathcal{O} &\rightarrow& H^0(A)\otimes
\mathcal{O}&\rightarrow
&H^0(B)\otimes \mathcal{O}&\rightarrow &0\\
&& \parallel && \downarrow&
&\downarrow& &\\
0&\rightarrow& W\otimes \mathcal{O} &\rightarrow& A&\rightarrow
& B&\rightarrow &0.\\
&& && \downarrow&
& \downarrow&&\\
&& && 0& & 0&&
\end{array}
\end{equation}

Actually, $B=M_{W,G}.$ Thus, we can state Mercat's results as
follows

\begin{theorem}\cite[Th\'{e}or\`{e}me B.1. and Th\'{e}or\`{e}me A.1]{me}\label{teoprin}
Let $G$ be a stable vector bundle of rank $s$ and degree $d_G>
2gs$. For a general generated subspace $W\subseteq H^0(G)$,
$M_{W,G}$ is stable and generated of slope $<2.$
\end{theorem}

From the construction of the stable bundles $ M_{W,G}$ we have the
following Proposition

\begin{proposition}\label{proppetri}
The Petri map for $(M_{W,G}, W^*)$ is injective.
\end{proposition}
\begin{proof}
The kernel of the Petri map
$$W^*\otimes H^0(M_{W,G}^*\otimes K)\rightarrow H^0(M_{W,G}\otimes
M_{W,G}^*\otimes K)$$ is $H^0(G^*\otimes M_{W,G}^*\otimes K).$

Since $M_{W,G}\otimes G$ is semistable and $\mu(M_{W,G}\otimes G)
>2g > 2g-2$,  $H^0(M_{W,G}^*\otimes G^*\otimes K)$. Thereby,
the Petri map is injective.
\end{proof}

The stable vector bundles $M_{W,G}$ in Theorem \ref{teoprin},
together with the stable bundles with $k\leq n$, and $0<d<2n$ were
tensored in \cite{bfmn}, by line bundles with sections to produce
Brill-Noether bundles with $d>2n$.

\begin{remark}\begin{em}\label{remmul}
Actually, let $d=d''+ nd'$ with $d'\geq (k'-1)(k'+g)/k'$,
$0<d''<2n$ and $1\leq k'\leq g$. It was proved in \cite[Theorem
3.2]{bfmn} that if $n\leq d'' + (n-k)g$ and $(d'',k)\not= (n,n),$
$B(n,d,k'k)\not= \emptyset .$ Furthermore, applying the Serre
duality, $B(n,2n(g-1)-d, k'k-d+n(g-1))\not= \emptyset .$ Moreover,
such bundles determine a region in the Brill-Noether map (see
\cite{bgn}) that extends beyond the region determined by the
stable bundles in \cite{mon1} (see \cite[Section 5]{bfmn}.
\end{em}\end{remark}

We want to study the coherent systems defined by the stable
bundles $M_{W,G}\otimes L$ with $L$ a line bundle with sections.

Let $(E,V)$ be a coherent system of type $(n,d,k)$ and $(F,W)$ a
subsystem of $(E,V)$. To study the stability of a coherent system
we can restrict to subsystems of the form $(F,H^0(F)\bigcap V)$
since $\mu _\alpha (F,W)\leq \mu _\alpha (F,H^0(F)\bigcap V)$ for
all $\alpha
>0$. Assume $E$ is stable and $k> n.$ If $h^0(F)\leq n_F$, $\mu
_{\alpha}(F,W)<\mu_{\alpha}(E,V), $ for all $\alpha
>0$. If $h^0(F)>n_F$ we have the following Lemma.

\begin{lemma}\label{prop1} Let $(F,W)$ be a coherent system
with $h^0(F)>n_F.$ If $H^0(F^*)\not= 0$ and $F\cong
\mathcal{O}^r\oplus G$ with $r\geq 1$ and $H^0(G^*)=0$,
$\frac{h^0(F)}{n_F}<\frac{h^0(G)}{n_G}$. Moreover, $\mu
_{\alpha}(F,W)<\mu_{\alpha}(G,H^0(G))$ for all $\alpha >0$.

\end{lemma}

\begin{proof}
If $F\cong \mathcal{O}^r\oplus G$ with $H^0(G^*)=0$, $r\geq 1$ and
$h^0(F)>n_F,$
\begin{equation}\label{equ}
\frac{h^0(F)}{n_F}=\frac{r+h^0(G)}{n_F}<\frac{h^0(G)}{n_F-r}=\frac{h^0(G)}{n_G}.
\end{equation}

Moreover, $d_F=d_G$ and $r+n_G=n_F<h^0(F)=r+h^0(G)$. Hence,
$\mu(F)<\mu(G)$ and from (\ref{equ}), $ \mu
_{\alpha}(F,W)<\mu_{\alpha}(G,H^0(G))$ for all $\alpha >0.$

\end{proof}

\begin{remark}\begin{em}\label{remstab011}
From Lemma \ref{prop1}, in order to prove the $\alpha$-stability of
 coherent system $(E,V)$ with $E$ stable we can just consider,
without loss of generality, subsystems $(F,V\bigcap H^0(F))$ that
satisfies $h^0(F)> n_F$, $H^0(F^*)=0.$
\end{em}\end{remark}

\section{Coherent systems of type $(n,d,n+1)$}

 For a general curve it was proved in \cite{lb} that
$G(\alpha:n,d,n+1)\not= \emptyset$ if and only if $g\geq
\beta(1,d,n+1) \geq 0$ and if $g=\beta$ , $n\nmid g$. Moreover, if
$\beta \geq 0$ then $G(n,d,n+1):= G(\alpha :n,d,n+1) = G(\alpha
':n,d,n+1)$ for all $\alpha ,\alpha ' >0.$ For $d> g+n$, the
geometry of the $G_L(n,d,n+1)$ was described in \cite{bun}.

For any curve we have the following Propositions. Note that there
is an overlap between them, but the proofs illustrate different
methods that could be used for $k>n+1.$

\begin{proposition}\label{teo11} If $d\geq n+g$ and $n\geq g$,
$U(n,d,n+1)\not= \emptyset .$
\end{proposition}
\begin{proof} Let $L$ be a line bundle of degree $d\geq 2g.$
Hence, $h^0(L)=d+1-g$. Moreover, it is generated and from
\cite[Proposition 3.2]{pr} a general subspace $V$ of dimension
$d+1-g\geq n+1\geq 2$ generates $L$. By Theorem \ref{teoprin},
$M_{V,L}$ is stable and from Remark \ref{remlineal1} and Remark
\ref{rem1} $(M_{V,L},V^*)$ is $\alpha$-stable for all $\alpha >0$.
Therefore, $U(n,d,n+1)\not= \emptyset .$
\end{proof}

\begin{proposition}\label{propn1}
If $n +g\leq d< 2n$ and $d'\geq 0$,
\begin{itemize}
\item $G(\alpha:n,d+nd' ,n+1)\not= \emptyset$ for all $\alpha >
0$; \item $U(n,d+nd',n+1)\not= \emptyset .$
\end{itemize}
\end{proposition}

\begin{proof}
For any $n+g\leq d< 2n$, there exists a generated stable bundle
$E$ (see Theorem \ref{teoprin}). Moreover, from Remarks
\ref{rem111} and \ref{rem11} there exists $V\subset H^0(E)$ such
that $(E,V)\in G_g(n,d,n+1)$ and it is $\alpha$-stable for all
$\alpha > 0$. Hence, $U(n,d,n+1)\not= \emptyset $. The Proposition
follows now from \cite[Lemma 1.5]{rag}.
\end{proof}

\begin{proposition}\label{proptensor1} For any $n\geq 2$ and
$d>2gn$, there exists $(L,W)\in G_g(1,d,n+1)$ such that $M_{W,L}$
is stable. Moreover, $U(n,d,n+1)\not= \emptyset $.
\end{proposition}
\begin{proof} Let $G$ be stable of rank $n$ and degree
$d>2gn.$ Hence, $G$ is generated and from \cite[Proposition
3.4]{pr} is generated by a subspace $W\subset H^0(G)$ of dimension
$n+1$.

The Proposition follows from the dual span correspondence since
the dual of the kernel of the evaluation map $M_{W,G}$ is a line
bundle $L$. Hence, $(G,W)$ is $\alpha _L$-stable since $(L,W)$ is
$\alpha$-stable for all $\alpha
>0$ (see \cite[Corollary 5.10]{bomn}). Moreover, $G$ is stable,
and hence $(G,W)$ is $\alpha $-stable for all $\alpha >0$.
Therefore, $U(n,d,n+1)\not= \emptyset .$
\end{proof}

\section{Coherent systems of type $(n,d,n+s)$}

To study $G(\alpha:n,d,n+s)$ with $n+sg +s'=d<2n$ for any $\alpha
>0$, we shall consider three cases, depending on $s'$, namely
\begin{enumerate}
\item $d< n+sg$;\item $d=n+sg$; \item $d>n+sg.$
\end{enumerate}

In this section we will give a complete description of the moduli
spaces $G(\alpha :n,d,k)$, for $d\leq n+sg$ (see Theorem
\ref{teo1} and \ref{teoempty}).

As we have seen (see Remark \ref{remstabmerc}), if $d<n+sg$,
$G_0(n,d,n+s)=\emptyset$ and the emptiness is related to the
non-existence of semistable bundles of type $(n,d,n+s)$. Our object
in this section is to generalize such relation to arbitrary $\alpha
>0$ and prove

\begin{theorem}\label{teoempty} If $d< n+gs,$
 $G(\alpha:n,d,n+s)= \emptyset$ for all $\alpha >0.$
\end{theorem}

\begin{remark}\begin{em}\label{remmio} For general curves, Theorem
\ref{teoempty} was proved in \cite[Corollary 3.10]{lb}.
\end{em}\end{remark}

We shall prove Theorems \ref{teoempty} and \ref{teo1} by means of a
sequence of Propositions.

The following Lemma follows at once from Proposition \ref{propm}
and Lemma \ref{prop1}.

\begin{lemma}\label{lematodos} Let $E$ and $F$ be vector bundles with
$\mu (F) \leq \mu (E) <2.$ If either $F$ is semistable or $E$
semistable and $F$ a subbundle of $E$,
\begin{equation}\label{eqtodos}
\frac{h^0(F)}{n_F}\leq \frac{\mu (E)-1}{g} +1.
\end{equation}
\end{lemma}

Recall from Remark \ref{remstab011} that if $E$ is stable, without
loss of generality, we can just consider subsystems $(F,V\bigcap
H^0(F))$ that satisfies $h^0(F)> n_F$, $H^0(F^*)=0.$

\begin{proposition}\label{propt1} Let $(E,V)$ be a coherent system of
type $(n,d,n+s)$ with $d\leq n+sg$,  $0<sg<n$. If $E$ is semistable,
$(E,V)$ is $\alpha$-semistable for all $\alpha >0$. Moreover, if $E$
is stable, $(E,V)$ is $\alpha$-stable for all $\alpha >0$.
\end{proposition}

\begin{proof} Let $(F,W)\subset (E,V)$ be a
subsystem with $H^0(F^*)=0$. Since $E$ is stable and $\mu (E) <2$,
from Lemma \ref{lematodos}
\begin{equation}\label{eq4}
\frac{\dim W}{n_F} \leq \frac{h^0(F)}{n_F}\leq 1+
\frac{\mu(E)-1}{g}\leq \frac{n+s}{n}.
\end{equation}

Therefore, from (\ref{eq4}) and the semistability of $E$, $\mu
_\alpha(F,W) \leq \mu_\alpha (E,V)$ for all $\alpha
>0$. If $E$ is stable $\mu
_\alpha(F,W) < \mu_\alpha (E,V)$ for all $\alpha
>0$.
\end{proof}

\begin{proposition}\label{propstability} If
$(E,V)$ is an $\alpha$-stable coherent system of type $(n,d,n+s)$
with $d\leq n+sg$, $E$ is stable.
\end{proposition}

\begin{proof}  Suppose $Q$ is a stable quotient bundle such
that $\mu(Q) \leq \mu(E)<2$.

From  Lemma \ref{lematodos}
$$\frac{h^0(Q)}{n_Q} \leq \frac{\mu(Q)-1}{g}+1 \leq
\frac{\mu(E)-1}{g}+1 \leq\frac{n+s}{n}.$$ Hence, $\mu _{\alpha
}(Q,W)\leq \mu _{\alpha}(E,V)$ which is a contradiction to the
$\alpha$-stability of $(E,V)$. Therefore, $E$ is stable.
\end{proof}

{\it Proof of Theorem \ref{teoempty}} If $(E,V)\in
G(\alpha:n,d,n+s)$, from Proposition \ref{propstability}, $E$ is
stable which a contradiction (see Remark \ref{remstabmerc}).
Therefore, $G(\alpha:n,d,n+s)= \emptyset$ for all $\alpha >0$.

$\hfill\Box$

For $d=n+sg$ we have the following Theorem.

\begin{theorem}\label{teo1}
Let  $2 \leq gs <n$ and $k=n+s$. If $d=n+gs$,
\begin{enumerate}
 \item $G(\alpha:n,d,k)\not= \emptyset$ for all
$\alpha >0$;
 \item $G(n,d,k): = G(\alpha: n, d, k) = G(\alpha':n,
d, k) $ for $\alpha ,\alpha ' >0$ i.e. $\alpha _L =0$;\item
$U(n,d,k)\not= \emptyset $. Moreover, $U(n,d,k)=G(n,d,k)\subset
G_g(n,d,k);$ \item $ G(n,d,k)$ is smooth irreducible of dimension
$ \beta $. Moreover, $G(n,d,k)\cong \mathcal{M}(s,d) $. \item
$G(\alpha :s,d,k)=G(n,d,k)$ for all $\alpha >0$
\end{enumerate}
\end{theorem}

\begin{proof}
 From Theorem  \ref{teoprin} there exists a
generated coherent system $(E,V)$ in $G_0(n,d,n+s) $ with $E$
stable and from Proposition \ref{propt1} $(E,V)$ is
$\alpha$-stable for all $\alpha >0$. Therefore,
$G(\alpha:n,d,n+s)\not= \emptyset $ for all $\alpha >0$. The
equality $G(\alpha; n, d, k) = G(\alpha';n, d, k) $ for $\alpha
,\alpha ' >0$ follows from Proposition \ref{propstability}. Part
$(4)$ follows from Proposition \ref{proppetri} and the dual span
correspondence gives the isomorphism $G(n,d,n+s)\cong
\mathcal{M}(s,d).$

Part $(5)$ follows from \cite[Corollary 5.10]{bomn} and the dual
span correspondence. Note that in this case, for the dual span
correspondence we do not need $(n,k)=1$, since the vector bundles
are stable and the case $\frac{k'}{n'}=\frac{k}{n}$ is allowed
(see \cite[Corollary 5.10]{bomn}).

\end{proof}

\section{Case $d>n+ sg$}

Assume now that $d=n+sg+s'$ with $0<s'<g$. From Remark
\ref{remdgrande}, $G_0(n,d,n+s)\not= \emptyset.$

For any $(E,V)\in G_0(n,d,n+s)$ we have the following Proposition.

\begin{proposition}\label{propdstab} Any $(E,V)\in G_0(n,d,n+s)$
 is generically generated.
\end{proposition}
\begin{proof} Let $(E,V)\in
G_0(n,d,n+s)$ and  $I_E$ be the image of the evaluation map
$V\otimes \mathcal{O}\rightarrow E$. If $n_I$ is the rank of
$I_E$, $n_I+r =n$ with $r\geq 0$. Let $I_E=\mathcal{O}^t\oplus N$
with $H^0(N^*)=0$ and $t\geq 0.$

Since, $E$ is semistable  $\mu(N)\leq \mu (E)$. From Lemma
\ref{prop1} and Proposition \ref{propm},
$$\frac{n+s}{n_I}\leq \frac{h^0(I_E)}{n_I}\leq \frac{h^0(N)}{n_N}$$
and
$$\frac{n_I+r+s}{n_I}\leq 1+ \frac{\mu (N)-1}{g} \leq 1+ \frac{\mu
(E)-1}{g}=1+\frac{s}{n}+\frac{s'}{ng}.$$ Hence,
$$\frac{r}{n_I}\leq \frac{s'}{ng} +\frac{s}{n}-\frac{s}{n_I} \leq
\frac{s'}{ng}. $$ Since $0<s'<g$, $0\leq r<1.$ Therefore, $r=0$
and hence $(E,V)$ is generically generated.
\end{proof}

Let $(B,H^0(B))\in G_0(n,d,n+s)$ be the coherent system defined in
section \S 3. Let $(F,H^0(F))\subset (B,H^0(B))$ be a coherent
subsystem. Without loss of generality (see Remark
\ref{remstab011}) we can assume that $h^0(F)> n_F$, $H^0(F^*)=0$.
Furthermore, the following Lemma allow us to assume also that $(F,
H^0(F))$ is generically generated.

\begin{lemma}\label{lemmasubgenerados} Let $(F,H^0(F))$ be a
subsystem of $(B,H^0(B))$ with $h^0(F)> n_F$ and  $H^0(F^*)=0$. If
$(F,H^0(F))$ is not generically generated then there exists a
generated subsystem $(N,H^0(N))\subseteq (B,H^0(B))$ with
$H^0(N^*)=0$ such $\mu _{\alpha}(F,H^0(F))<\mu_{\alpha}(N,H^0(N))$
for all $\alpha
>0$.
\end{lemma}

\begin{proof} If $h^0(F)=n_F+b$, from Proposition \ref{propm}, $n_F+bg\leq d_F$
i.e. $d_F= n_F+bg+b'$ with $g>b'\geq 0$.

Assume $(F,H^0(F))$ is not generically generated and let $I_F$ be
 the image of the evaluation map $H^0(F)\otimes
\mathcal{O}\rightarrow F$. If $n_I$ is the rank of $I_F$, $n_I+r
=n_F$ with $r\geq 1$.

Let $I_F=\mathcal{O}^t\oplus N$ with $H^0(N^*)=0$ and $t\geq 0.$
Note that $(N,H^0(N))$ is generated.

As in the proof of Proposition \ref{propdstab}, from Lemma
\ref{prop1} we have that,
\begin{equation}\label{eqsub00}\frac{h^0(F)}{n_F}<\frac{n_F+b}{n_I}= \frac{h^0(I_F)}{n_I}\leq
\frac{h^0(N)}{n_N}.\end{equation} Suppose $\mu(N)\leq \mu (F)$.
From Proposition \ref{propm},
$$\frac{n_I+r+b}{n_I}\leq 1+ \frac{\mu (N)-1}{g} \leq 1+ \frac{\mu
(F)-1}{g}=1+\frac{b}{n_F}+\frac{b'}{n_Fg}.$$ Hence,
$$\frac{r}{n_I}<\frac{b'}{n_Fg} +\frac{b}{n_F}-\frac{b}{n_I} <
\frac{b'}{n_Fg}.$$

Since $b'<g$, $r$ must have to be $<1$, which is a contradiction.
Therefore, $\mu(F)<\mu (N)$. This last inequality together with
$(\ref{eqsub00})$ implies that $\mu
_{\alpha}(F,H^0(F))<\mu_{\alpha}(N,H^0(N))$ for all $\alpha
>0$.
\end{proof}

Therefore, from Lemmas \ref{prop1} and \ref{lemmasubgenerados}, to
prove the $\alpha$-stability of $(B,H^0(B))$ we can assume that
the subsystems $(F,H^0(F))$ are generically generated with
$h^0(F)>n$ and $H^0(F^*)=0$.

\begin{proposition}\label{propdg} $(B,H^0(B))$ is $\alpha $-stable for all $\alpha
>0$. Moreover, $U(n,d,n+s)\not= \emptyset .$
\end{proposition}

\begin{proof} Let $(F,H^0(F))$ be a generically generated subsystem of
type $(n_F, n_F+bg+b',n_F+b)$ with $H^0(F^*)=0$. Since $B$ is
stable, to prove the $\alpha$-stability of $(B,H^0(B))$ we need to
prove that $\frac{b}{n_F}\leq \frac{s}{n}.$

 From the stability of $B$
\begin{equation}\label{eqslope}
1+\frac{bg}{n_F} + \frac{b'}{n_F} =\mu(F)< \mu(B)=1+\frac{sg}{n} +
\frac{s'}{n}.
\end{equation}
Thus, if $\frac{s'}{n}\leq \frac{b'}{n_F}$, from
$(\ref{eqslope})$, $\frac{b}{n_F}<  \frac{s}{n}$ and hence, $\mu
_{\alpha} (F,H^0(F))< \mu _{\alpha}(B,H^0(B))$ for all $\alpha
>0.$

Assume $\frac{s'}{n}> \frac{b'}{n_F}$.

Since $(I_F,H^0(F))$ is generated, the exact sequence
$(\ref{eqsystems1})$ fits in the following diagram
\begin{equation}\label{diagram2g}
\begin{array}{ccccccccc}
&& 0&& 0&
& 0&&\\
&&\downarrow && \downarrow&
& \downarrow&&\\
0&\rightarrow& K_F &\rightarrow& H^0(F)\otimes
\mathcal{O}&\rightarrow
&I_F&\rightarrow &0\\
&& \downarrow && \downarrow&
&\downarrow& &\\
0&\rightarrow& M_B^* &\rightarrow&H^0(B)\otimes
\mathcal{O}&\rightarrow
& B&\rightarrow &0.\\
\end{array}
\end{equation}

By construction, $M_B$ is a stable bundle of degree $d>2gs$ and
rank $s$. From the stability of $M_B$ and diagram
$(\ref{diagram2g})$ we have
$$\frac{n+sg+s'}{s}=\mu (M_B)< \mu
(K^*_F)=\frac{d_{{I_F}}}{b}<\frac{n_F+bg+b'}{b}.
$$
That is,
\begin{equation}\label{eqslopedes}\frac{n+s'}{s}<\frac{n_F+b'}{b}.
\end{equation}

If \begin{equation}\label{eqeq}\frac{b}{n_F}> \frac{s}{n},
\end{equation}
$s'>b'\frac{n}{n_F}>\frac{b's}{b} \, $ i.e. $\frac{s'}{s}>
\frac{b'}{b}.$ But, from $(\ref{eqslopedes}),$
\begin{equation}\label{eqslopedesi}
\begin{array}{lll}
\frac{n}{s}&<& \frac{n_F}{b} +\frac{b'}{b} - \frac{s'}{s}\\
&&\\
 &<&\frac{n_F}{b},
\end{array}
\end{equation}
which is a contradiction to $(\ref{eqeq})$. Hence,
$\frac{b}{n_F}\leq \frac{s}{n}$ and, from the stability of $B$,
$\mu _{\alpha} (F,H^0(F))< \mu _{\alpha}(B,H^0(B))$ for all
$\alpha
>0.$ Therefore, $(B,H^0(B))$ is $\alpha $-stable  for all $\alpha
>0$ and $U(n,d,n+s)\not= \emptyset .$
\end{proof}

\begin{theorem}\label{teo2} If
 $k=n+s$ and $d=n+sg+s'$ with $0<s'<g$ and $2 \leq gs <n$,
\begin{enumerate} \item $G(\alpha:n,d,k)\not= \emptyset$ for all
$\alpha >0$; \item $U(n,d,k)\not= \emptyset$. Moreover, any
$(E,V)\in U(n,d,k)$ is generically generated; \item $ G(\alpha
:n,d,k)$ has a smooth irreducible component $G^0(n,d,k)$ of
dimension $ \beta $. Moreover, $G^0(n,d,k)$ is birationally
equivalent to the Grassmanian bundle $Grass (s')$ over
$\mathcal{M}(s,d)$.
\end{enumerate}
\end{theorem}

\begin{proof} Part (1) and (2)  follows from Propositions
\ref{propdg} and \ref{propdstab}, respectively. The first part of
$(3)$ follows from Proposition \ref{proppetri}.

For any $G\in \mathcal{M}(s,d)$, let $U_G$ be the set
$$U_G=\{V\in Grass(s',H^0(G)) : (G,V)\,
 {\mbox{ is generated with }} \,  M_{V,G} {\mbox{ stable}} \}.$$

By openness of stability, $U_G$ is an open set of
$Grass(s',H^0(G))$. The open sets $U_G$ for all $G\in
\mathcal{M}(n,d)$ define an open set $Z$ in the Grassmannian
bundle ${\it Grass}(s)$ over $\mathcal{M}(n,d)$. The dual span
correspondence define a coherent system in $U(n,d,n+s)$. From the
universal properties of the moduli space $G_L(n,d,k),$ the map
from $Z$ to $G^0(n,d,k)$ is regular and hence it gives a
birrational equivalence.
\end{proof}

From Theorem \ref{teo1} and \ref{teo2} we conclude:

\begin{corollary}\label{corprincipal}
For any Brill-Noether triple $(n,d,k)$ with $n<d<2n$ and $k>n$
with $n\leq d+(n-k)g$, $U(n,d,k)\not= \emptyset$. Moreover, there
is an open set $Z$ of $B(n,d,k)$ such that any $E\in Z$ defines a
$\alpha$-stable coherent system  for all $\alpha
>0$.
\end{corollary}

\section{ Tensoring coherent systems}

Given two coherent systems $(E_i,V_i)$ of type $(n_i,d_i,k_i)$,
with $i=1,2$, the pair $(E_1\otimes E_2, V_1\otimes V_2)$ needs
not to be a coherent system of type
$(n_1n_2,d_1n_2+d_2n_1,k_1k_2)$.

However, to get a coherent system of type
$(n_1n_2,d_1n_2+d_2n_1,k_1k_2)$ tensor the associated sequences
$(\ref{eqsystems1})$ and $(\ref{eqsystems2})$ of $(E_i,V_i)$ with
$E_j$ to get the following exact sequences
\begin{equation}\label{eqsystems11}
0\rightarrow K_{{E_i}}\otimes E_j \rightarrow V_i\otimes
E_j\rightarrow I_{{E_i}}\otimes E_j \rightarrow 0
\end{equation}
and
\begin{equation}\label{eqsystems22}
0\rightarrow I_{{E_i}}\otimes E_j \rightarrow E_i\otimes
E_j\rightarrow H\otimes E_j\oplus \tau  \rightarrow 0.
\end{equation}

If $H^0(K_{{E_i}}\otimes E_j)=0$, $$ (\varphi _{(E_i,V_i)}\otimes
id)^*: V_i\otimes H^0(E_j)\rightarrow H^0(E_i\otimes E_j)$$ is
injective  and hence we identify $V_i\otimes V_j$ with the image
of $V_i\otimes V_j$ in $H^0(E_i\otimes E_j)$ under the
 inclusions $V_i\otimes V_j \hookrightarrow V_i\otimes H^0(E_j)
 \hookrightarrow H^0(E_i\otimes E_j)$.  In this case, we define the {\it tensor
product} of $(E_j,V_j)$ by $(E_i,V_i)$ to be the coherent system
$(E_i\otimes E_j,V_i\otimes V_j)$ of type
$(n_1n_2,d_1n_2+d_2n_1,k_1k_2)$.

In this section we are interested in tensor coherent systems
$(E,V)$ of type $(n,d,k)$ by coherent systems $(L,W)$ of type
$(1,d',k').$

Given $(L,W)$ of type $(1,d',k')$ we have the exact sequences
\begin{equation}\label{eqim01} 0 \rightarrow K_L \rightarrow
W\otimes \mathcal{O} {\rightarrow} I_L\rightarrow 0
\end{equation}
and
$$0\rightarrow I_L \rightarrow L {\rightarrow} \tau \rightarrow
0,$$ where $I_L$ is the image of the evaluation map $W\otimes
\mathcal{O}\rightarrow L$.

 Assume $H^0(K_L\otimes E)=0$ and let
$(E\otimes L, V\otimes W)$ be a coherent system of type
$(n,d+nd'+kk')$. To study the stability of $(E\otimes L , V\otimes
W)$ we want to describe $H^0(F)\cap (V\otimes W)$ for any
subbundle $F\subseteq E\otimes L$. The following Lemmas will do
it.

First recall that any subbundle $F\subseteq E\otimes L$ defines a
subbundle $F' \subseteq E\otimes I_L$ that fits in the following
diagram
\begin{equation}\label{diagimag}
\begin{array}{ccccccccc}
&& 0 && 0&
&0& &\\
&& \downarrow&& \downarrow&
&\downarrow& &\\
0&\rightarrow& F' &\rightarrow& F&{\rightarrow}
&\tau '&\rightarrow &0\\
&& \downarrow&& \downarrow&
&\downarrow& &\\
0&\rightarrow& I_L\otimes E &\rightarrow& L\otimes E&{\rightarrow}
&\tau &\rightarrow &0.\\
\end{array}
\end{equation}

\begin{lemma}\label{lemasubi} $H^0(F)\cap (V\otimes W)= H^0(F')\cap (V\otimes W).$
\end{lemma}
\begin{proof} From the diagram $(\ref{diagimag})$,
$H^0(F')=H^0(F)\cap H^0(I_L\otimes E)$. Hence,
$$
\begin{array}{ccl}
H^0(F')\cap (V\otimes W)&=&H^0(F)\cap H^0(I_L\otimes E)\cap
(V\otimes W)\\
 &=&H^0(F)\cap (V\otimes W),
\end{array}
$$
since $H^0(K_{{L}}\otimes E)=0$ and $V\otimes W \subseteq
H^0(E)\otimes W \subseteq H^0(I_L\otimes E).$
\end{proof}

 For the subbundle $F'\subset E\otimes I_L$, we have
the exact sequences $ 0\rightarrow F' \stackrel{i}{\rightarrow}
I_L\otimes E \stackrel{q}{\rightarrow} Q\rightarrow 0$ and
\begin{equation}\label{eqlineal2}
0\rightarrow F'\otimes I_L^* \rightarrow E
\stackrel{p}{\rightarrow} Q\otimes I_L^*\rightarrow 0
\end{equation}
where $Q$ is the quotient bundle.

Tensor $(\ref{eqlineal2})$ by the sequence $(\ref{eqim01})$ to get
the following diagram
\begin{equation}\label{diag1}
\begin{array}{ccccccccc}
&& 0 && 0&
&0& &\\
&& \downarrow&& \downarrow&
&\downarrow& &\\
0&\rightarrow& K_L\otimes F'\otimes I_L^* &\rightarrow& K_L\otimes
E&\stackrel{id\otimes p}{\rightarrow}
&K_L\otimes Q\otimes I_L^*&\rightarrow &0\\
&& \downarrow&& \downarrow&
&\downarrow& &\\
0&\rightarrow& W\otimes F'\otimes I_L^* &\rightarrow& W\otimes
E&\stackrel{id \otimes p}{\rightarrow}
&W\otimes Q\otimes I_L^*&\rightarrow &0\\
&& \downarrow&& \downarrow&
&\downarrow& &\\
 0&\rightarrow& F' &\rightarrow& I_L\otimes E&\stackrel{q}{\rightarrow}
&Q&\rightarrow &0\\
&& \downarrow&& \downarrow&
&\downarrow& &\\
&& 0 && 0&
&0& &\\
\end{array}
\end{equation}

 From the cohomology  of
(\ref{diag1}), we have the following diagram
\begin{equation}\label{diag2}
\begin{array}{ccccccccc}
&& && &
&0& &\\
&& && &
&\downarrow& &\\
&& 0&& 0&\rightarrow
&H^0(K_L\otimes Q\otimes I_L^*)&\rightarrow &\\
&& \downarrow&& \downarrow&
&\downarrow& &\\
0&\rightarrow& W\otimes H^0(F'\otimes I_L^*) &\rightarrow&
W\otimes H^0(E)&\stackrel{id \otimes p^*}{\rightarrow}
&W\otimes H^0(Q\otimes I_L^*)& \rightarrow&\\
&& \downarrow&& \downarrow&
&\downarrow& &\\
 0&\rightarrow& H^0(F') &\stackrel{i^*}{\rightarrow}& H^0(I_L\otimes
E)&\stackrel{q}{\rightarrow}
&H^0(Q)&\rightarrow &\\
&& \downarrow&& \downarrow&
&\downarrow& &\\
\end{array}
\end{equation}
since $H^0(K_L\otimes E)=0$.

\begin{lemma}\label{propmulti} $ H^0(F')\cap (W\otimes H^0(E))=W\otimes
H^0(F'\otimes I_L^*).$
\end{lemma}

\begin{proof}
From the commutativity of diagram (\ref{diag2}), $$W\otimes
H^0(F'\otimes I_L^*)\subseteq (W\otimes H^0(E))\cap H^0(F').$$ To
prove $(W\otimes H^0(E))\cap H^0(F')\subseteq W\otimes
H^0(F'\otimes I_L^*)$ take any $0\not=a\in W$. Hence, $<a>\otimes
H^0(E) \subset W\otimes H^0(E)$. From $(\ref{diag2})$ we have the
following diagram

\begin{equation}\label{diag3}
\begin{array}{ccccccccc}
&& 0&& 0&
&& &\\
&& \downarrow&& \downarrow&
&& &\\
0&\rightarrow& <a>\otimes H^0(F'\otimes I_L^*) &\rightarrow&
<a>\otimes H^0(E)&\stackrel{id \otimes p^* }{\rightarrow}
&<a>\otimes H^0(Q\otimes I_L^*)&\rightarrow &\\
&& \downarrow  && \downarrow \gamma&
&\downarrow \beta & &\\
 0&\rightarrow& H^0(F') &\stackrel{i^*}{\rightarrow}& H^0(I_L\otimes E
)&\stackrel{q^*}{\rightarrow}
&H^0(Q)&\rightarrow &\\
\end{array}
\end{equation}

Let $\{s_i\} $ be a basis of $H^0(E)$. Let
 $$c= \lambda a\otimes \sum \nu _i s_i \in \, <a>\otimes
H^0(E)$$ be such that $ \gamma (c) \in i^* (H^0(F'))\subset
H^0(I_L\otimes E).$
 From diagram (\ref{diag3}),
$q^*(\gamma (c))=0 \in H^0(Q)$ i.e. $q^*(\gamma (c))$ is the zero
section. Hence,
$$q^*(\gamma (c))=\beta(id \otimes p^*)(c)= \beta (\lambda a\otimes
p^*(\sum \nu _i s_i))=0.$$ That is, for
 $x\in X$
\begin{equation}
x \mapsto \beta(id \otimes p^*)(c)(x)=\lambda a(x)\otimes p^*(\sum
\nu _i s_i)(x)=0.
\end{equation}

Since $a\in W$ is a non-zero section, $\lambda a(x)=0$ only for a
finite set of points $x_j\in X$. Hence, $p^*(\sum \nu _i s_i)$
must be the zero section. That is, $\sum \nu _i s_i \in \ker \, (
p^*)= H^0(F'\otimes I_L^*)$. Therefore, if $ \gamma (c) \in i^*
(H^0(F'))$, $c\in <a>\otimes H^0(F'\otimes I_L^*)$. Thereby, for
all $a \in W$,
$$(<a>\otimes
H^0(E))\cap H^0(F') \subseteq <a>\otimes H^0(F'\otimes
I_L^*)\subseteq W\otimes H^0(F'\otimes I_L^*)$$ and hence,
$$(W\otimes H^0(E))\cap H^0(F')\subseteq W\otimes H^0(F'\otimes I_L^*).$$
Therefore, $W\otimes H^0(F'\otimes I_L^*)=(W\otimes H^0(E))\cap
H^0(F').$
\end{proof}

Moreover, for $V\varsubsetneq H^0(E)$ we have the following Lemma

\begin{lemma}\label{propmulti2}
$H^0(F)\cap (W\otimes V)= W\otimes (H^0(F'\otimes I_L^*)\cap V).$
\end{lemma}

\begin{proof} From the injectivity of the first two columns in
diagram \ref{diag2} it is clear that $$W\otimes (H^0(F'\otimes
I_L^*)\cap V)\subseteq H^0(F')\cap(V\otimes W).$$

Moreover, from Lemma \ref{propmulti}
\begin{equation}
\begin{array}{lll}
 H^0(F')\cap(V\otimes W)&=& H^0(F')\cap ((H^0(E)\otimes
 W)\cap (V\otimes W))\\
 &=&(H^0(F')\cap (H^0(E)\otimes
 W))\cap (V\otimes W)\\
 &=&(H^0(F'\otimes I_L^*)\otimes W) \cap (V\otimes W)\\
 &\subseteq &(H^0(F'\otimes I_L^*)\cap V)\otimes W.
 \end{array}
 \end{equation}
 Therefore, $H^0(F')\cap (V\otimes W)= W\otimes (H^0(F'\otimes I_L^*)\cap V)$
and the Proposition follows from Lemma \ref{lemasubi}.
\end{proof}

For coherent systems $(L,\mathbb{C})$ of type $(1,d',1)$ we know
(see \cite[Lemma 1.5]{rag}) that $(E,V)$ is $\alpha $-stable if
and only if $(E\otimes L, V)$ is $\alpha $-stable. For coherent
systems of type $(1,d',k')$ we have the following Lemma.

\begin{lemma}\label{proplemma} Let $(E,V)$ be a coherent system of
type $(n,d,k)$ and $(L,W)$ of type $(1,d',k') $. Assume
$H^0(K_L\otimes E)=0$. Then $(E,V)$ is $k'\alpha$-stable if and
only if $(E\otimes L,V\otimes W)$ is $\alpha$-stable.
 Moreover, \begin{enumerate} \item if $(E,V)$ is $\alpha _L$-stable,
  $(E\otimes L,V\otimes W)$ is $\alpha _L$-stable;
  \item if $ U(n,d,k)\not= \emptyset$, $U(n,d+nd', kk')\not= \emptyset.$
  \end{enumerate}
\end{lemma}
\begin{proof} Note that if $(F,
H^0(F)\cap (V\otimes W))$ is a subsystem of $(E\otimes L,V\otimes
W)$, $(F'\otimes I_L^*, H^0(F'\otimes I_L^*)\cap V)$ is a coherent
subsystem of $(E,V)$, where $F'$ is the subbundle of $I_L\otimes
E$ that fits in $(\ref{diagimag})$.

The first part of the Lemma follows at once from the definitions
and Lemmas \ref{propmulti} and \ref{propmulti2} once we notice
that $\mu(F) -\mu(E\otimes L) \leq \mu(F'\otimes I_L^*)- \mu(E)$
and that for any subbundle $F\subseteq E\otimes L$
\begin{equation}\label{eqdimen} \frac{\dim ((W\otimes V)\cap
H^0(F))}{n_F}=\frac{k'\dim (H^0(F'\otimes I_L^*)\cap V)}{n_F}.
\end{equation}

 For the second part, if $(E,V)$ is $\alpha
_L$-stable,
\begin{equation}\label{eqcasifin}
0\leq \frac{k}{n}-\frac{\dim (H^0(F'\otimes I_L^*)\cap V)}{n_F}.
\end{equation}

From Lemmas \ref{propmulti} and \ref{propmulti2}
\begin{equation}\label{eqfinc}
\begin{array}{lll}
\mu (F)-\mu (E\otimes L)&\leq &\mu(F'\otimes I_L^*)-\mu(E)\\
&&\\
&<& \alpha (\frac{k}{n}-\frac{\dim (H^0(F'\otimes I_L^*)\cap V)}{n_F})\\
&&\\
&<& \alpha k'(\frac{k}{n}-\frac{\dim (H^0(F'\otimes I_L^*)\cap V)}{n_F})\\
&&\\
 &=&\alpha (\frac{k'k}{n}-\frac{\dim ((H^0(F'\otimes I_L^*)\cap V)\otimes
W)}{n_F})\\
&&\\ &=&\alpha (\frac{k'k}{n}-\frac{\dim (H^0(F)\cap (V\otimes
W))}{n_F}).
\end{array}
\end{equation}

Therefore, $\mu _{\alpha}(F,H^0(F)\cap (V\otimes W)) <\mu
_{\alpha}(E\otimes L, V\otimes W)$ for $\alpha >\alpha _L$ and
hence $(E\otimes L,V\otimes W)$ is $\alpha $-stable for all
$\alpha >\alpha _L$. It is well known that if $E$ is stable,
$E\otimes L$ is stable. Thereby, if $(E,V)\in U(n,d,k)$,
$(E\otimes L, V\otimes W)\in U(n,d+nd', kk').$

\end{proof}
\smallskip

\begin{theorem}\label{teodetodos} Suppose $U(n,d,k)\not=
\emptyset$. If there exists a coherent system of type $(1,d',k')$
with $K_L$ semistable and $(k'-1)d<nd'$, $U(n,d+nd', k'k)\not=
\emptyset .$
\end{theorem}
\begin{proof} Let $(E,V)\in U(n,d,k)$. By hypothesis, $K_L\otimes
E$ is semistable of slope $<0$, so $H^0(K_L\otimes E)=0$. The
Theorem follows at once from Lemma \ref{proplemma}.
\end{proof}

\begin{corollary}\label{teogrande2} Let $d'\geq 2g$, $k'\geq
2$  with $d(k'-1)<d'n$. If $U(n,d,k)\not= \emptyset$, $U(n,d
+nd',k'k)\not= \emptyset .$
\end{corollary}

\begin{proof}  Let $L$ be a line bundle of degree $d'\geq 2g$ and
$W$ a general subspace of $H^0(L)$ of dimension $k'\geq 2$ that
generates $L$. From Theorem \ref{teoprin}, $K_L=M^*_{W,L}$ is
stable and has slope $\frac{d'}{k'-1}.$ Hence, if $(E,V)\in
U(n,d,k)$, $H^0(K_L\otimes E)=0$, therefore,
 $(E\otimes L, V\otimes W)\in U(n,d
+nd',k'k)\not= \emptyset.$
 \end{proof}

Assume now that  $0<d''<2n$ and $k''\geq 1$ with $n\leq
d''+(n-k'')g$. If $(n,d'',k'')\not= (n,n,n)$, $U(n,d'',k'')\not=
\emptyset $ (see Corollary \ref{corprincipal} and \cite[Theorem
A]{bommn}). Hence, from Theorem \ref{teodetodos}, we have the
following Corollary.

\begin{corollary}\label{teogrande3} Suppose $0<d''<2n$ and $k''\geq 1$
with $n\leq d''+(n-k'')g$. Let $d'\geq 2g$, $k'\geq 2$  with
$d''(k'-1)<d'n$. If $(n,d'',k'')\not= (n,n,n)$, $U(n,d''
+nd',k'k'')\not= \emptyset .$
\end{corollary}

In the definition of the tensor product of coherent systems and in
the proof of Lemma \ref{propmulti} we use that $H^0(K_L\otimes
E)=0$. Such a condition seems very strong; however, for coherent
systems $(E,V)$ with $V=H^0(E)$ we have the following Proposition.

\begin{proposition}\label{proptensor01} Let $(E,V)$ be a
coherent system with $V=H^0(E)$ and $(L,W)$ a coherent system of
type $(1,d',k')$. \begin{enumerate} \item If $(E,V)$ is generated
and $H^0(L\otimes M^*_{E})=0$, $H^0(K_L\otimes E)=0.$ \item If
$(E,V)$ is injective, $H^0(K_L\otimes E)=0.$
\end{enumerate}
\end{proposition}
\begin{proof} Suppose $(E,V)$ is generated. Tensor the exact sequence
\begin{equation}
0\rightarrow M^*_E\rightarrow H^0(E)\otimes \mathcal{O}
\rightarrow E \rightarrow 0
\end{equation}
with
\begin{equation}
0\rightarrow K_L\rightarrow W\otimes \mathcal{O} \rightarrow I_L
\rightarrow 0
\end{equation}
to get the following diagram
\begin{equation}\label{diagcero}
\begin{array}{ccccccccc}
&& 0 && 0&
&0& &\\
&& \downarrow&& \downarrow&
&\downarrow& &\\
0&\rightarrow& K_L\otimes M^*_{E} &\rightarrow& W\otimes M^*_{E}
&{\rightarrow}
&I_L\otimes M^*_E&\rightarrow &0\\
&& \downarrow&& \downarrow&
&\downarrow& &\\
0&\rightarrow& K_L\otimes H^0(E) &\rightarrow& W\otimes H^0(E)
&{\rightarrow}
&I_L\otimes H^0(E)&\rightarrow &0\\
&& \downarrow&& \downarrow&
&\downarrow& &\\
 0&\rightarrow& K_L\otimes E &\rightarrow& W\otimes E&{\rightarrow}
&I_L\otimes E&\rightarrow &0.\\
&& \downarrow&& \downarrow&
&\downarrow& &\\
&& 0 && 0&
&0& &\\
\end{array}
\end{equation}
The first part of the Proposition follows from the cohomology of
diagram $(\ref{diagcero})$ once we notice that
$H^0(K_L)=0=H^0(M^*_E)$ and if $H^0(L\otimes M^*_{E})=0$,
$H^0(I_L\otimes M^*_{E})=0$.

 If $(E,V)$ is an injective, the Proposition  follows at once
 from the cohomology of the
following diagram
\begin{equation}\label{diagcero1}
\begin{array}{ccccccccc}
&& 0 && 0&
&0& &\\
&& \downarrow&& \downarrow&
&\downarrow& &\\
0&\rightarrow& K_L\otimes H^0(E) &\rightarrow& W\otimes H^0(E)
&{\rightarrow}
&I_L\otimes H^0(E)&\rightarrow &0\\
&& \downarrow&& \downarrow&
&\downarrow& &\\
 0&\rightarrow& K_L\otimes E &\rightarrow& W\otimes E&{\rightarrow}
&I_L\otimes E&\rightarrow &0\\
&& \downarrow&& \downarrow&
&\downarrow& &\\
&& 0 && 0&
&0& &\\
\end{array}
\end{equation}
since $H^0(K_L)=0$.
\end{proof}

Recall from \cite[Lemma 2.6]{lau} that if $B(n,d,k) \not=
\emptyset $ there exists $E\in B(n,d,k)$ with $h^0(E)=k$. Hence,
if $G_0(n,d,k)\not= \emptyset $, there exists
 $(E,V)\in G_0(n,d,k)$ with $V=H^0(E)$.

The following Theorems follow from Lemma \ref{proplemma} and
Theorem \ref{teodetodos}. The hypothesis in each one allow us to
define the tensor products of coherent systems of type $(n,d,k)$
by those of type $(1,d',k')$.

\begin{theorem}\label{teogrande02} Assume $0<d''<2n$ and $k''\geq 1$
with $n\leq d''+(n-k'')g$ and $(n,d'',k'')\not= (n,n,n)$. Let
$d=d''+nd'$ and $k=k'k''$ with $d'<2g$ and $1\leq k'$. If $\beta
(1,d',k')>\geq 0$, $U(n,d,k)\not= \emptyset .$
\end{theorem}

\begin{proof} If $\beta
(1,d',k')\geq 0$, there exist coherent systems $(L,W)$ of type
$(1,d',k')$.

For $k\leq n$, from \cite[Theorem A]{bommn}, $ U(n,d'',k'')\not=
\emptyset $ and any $(E,V)\in U(n,d'',k'')$ is injective. From
\cite[Lemma 2.6]{lau} there exists $(E,V)\in U(n,d'',k'')$ with
$V=H^0(E)$ and from Proposition \ref{proptensor01},
$H^0(K_L\otimes E)=0.$ From Lemma \ref{proplemma}, $U(n,d,k)\not=
\emptyset .$

For $k>n$, from the hypothesis and the proof of Theorems
\ref{teo1} and \ref{teo2} there exist generated coherent systems
$(C,V)\in U(n,d'',k'')$ with $V=H^0(C)$ and $M_C$ stable of slope
$>2g.$ Hence, for any line bundle $L$ of degree $d'\leq 2g$,
$H^0(M^*_C\otimes L)=0$ and from Proposition \ref{proptensor01},
$H^0(K_L\otimes E)=0$. Therefore, from Lemma \ref{proplemma}
$U(n,d,k)\not= \emptyset .$
\end{proof}

\begin{corollary}\label{corfinal} Assume $0<d''<2n$ and $k''\geq 1$
with $n\leq d''+(n-k'')g$ and $(n,d'',k'')\not= (n,n,n)$. Let
$d=d''+nd'$ and $k=k'k''$ with $0<d'<2g$ and $1\leq k'$. If $\beta
(1,d',k')\geq 0$, for any Brill-Noether triple $(n,d,k)$,
$U(n,d,k)\not= \emptyset$. Moreover, there is an open set $Z$ of
$B(n,d,k)$ such that any $E\in Z$ defines a $\alpha$-stable
coherent system for all $\alpha >0$.
\end{corollary}

 For hyperelliptic curve we have the following Theorem.

\begin{theorem}\label{teohyp}  Let $X$
be a hyperelliptic curve of genus $g\geq 3$. Assume $0<d''<2n$ and
$k''\geq 1$ with $n\leq d''+(n-k'')g$ and $(n,d'',k'')\not=
(n,n,n)$. If $d'=2(k'-1)$, $U(n,d+nd',k'k)\not=\emptyset .$
\end{theorem}
\begin{proof} Let $L$
be the hyperelliptic line bundle over $X$. For $1\leq k'\leq g$,
$h^0(L^{\otimes (k'-1)})=k'$ and the degree of $L^{\otimes
(k'-1)}$ is $d'=2(k'-1)$. As in Theorem \ref{teogrande02} we get
$U(n,d+nd',k'k)\not=\emptyset $ after tensoring with
$(L,H^0(L^{\otimes (k'-1)})).$
\end{proof}
\bigskip

 {\small \noindent \textbf{Acknowledgements.} The first author thanks
 P. Del Angel and F. Ongay  for useful
conversations and acknowledges the support of CONACYT grant
48263-F. The second author thanks CIMAT, where part of the work
was carried out, for their hospitality and acknowledges the
support of CONACYT grant 61011. The first author is a member of
the research group VBAC (Vector Bundles on Algebraic Curves).}


\end{document}